\documentclass[12pt]{article}
\usepackage{latexsym}
\usepackage{amsmath}
\usepackage{amssymb}
\usepackage{amsfonts}

\usepackage{url}

\def\claim#1{\begin{trivlist}\item[\hskip\labelsep\bf#1]\it}
\def\endclaim{\end{trivlist}}

\numberwithin{equation}{section}

\newtheorem{theorem}{Theorem}[section]

\newtheorem{pr}{Proposition}[section]

\newtheorem{remark}[theorem]{Remark}

\newcommand{\eproof}{{\mbox{\ }~\hfill
\mbox{\large $\Box$} \par \vskip 10pt}}

\title{Unique continuation property of solutions to general second order elliptic systems}
\author{N. Honda\thanks{Department of Mathematics, Hokkaido University, Sapporo 060-0810, Japan. (Email:
honda@math.sci.hokudai.ac.jp)}\qquad C.-L. Lin\thanks{Department of Mathematics, NCTS, National Cheng-
Kung University, Tainan 701, Taiwan.  Partially supported by the
Ministry of Science and Technology of Taiwan. (Email:
cllin2@mail.ncku.edu.tw)}\qquad G. Nakamura\thanks{Department of
Mathematics, Inha University, Incheon 402-751, Republic of Korea.
Partially supported by Korea Research Foundation.
(Email: 213028@inha.ac.kr)}\qquad S. Sasayama\thanks{Department of
Mathematics, Hokkaido University, Sapporo 060-0810, Japan. (Email:
sasayama@math.sci.hokudai.ac.jp)}\\
\\
Dedicated to the memory of Prof. Kenjiro Okubo}

\begin{document}
\date{}
\renewcommand{\theequation}{\thesection.\arabic{equation}}
 \maketitle
\begin{abstract}
This paper concerns about the weak unique continuation property of solutions of a general system of differential equation/inequality with a second order strongly elliptic system as its leading part. We put not only some natural assumption which we call {\sl basic assumptions}, but also some technical assumptions which we call {\sl further assumptions}. It is shown as usual by first applying the Holmgren transform to this inequality and then establishing a Carleman estimate for the leading part of the transformed inequality. The Carleman estimate given via a partition of unity and Carleman estimate for the operator with constant coefficients obtained by freezing the coefficients of the transformed leading part at a point. A little more details about this are as follows. Factorize this  operator with constant coefficients into two first order differential operators. Conjugate each factor by a Carleman weight and derive an estimate which is uniform with respect to the point at
 which we froze the coefficients for each conjugated factor by constructing a parametrix for its adjoint operator.
\end{abstract}

\section{Introduction}\label{sec1}
\setcounter{equation}{0}

Let $N\in\mathbb{N}$ and $\Omega$ be a domain in $\mathbb{R}^n$ with $n\geq 2$.
Suppose a ${\mathbb{C}}^n$ valued vector function $u(x)=(u_1(x), \cdots, u_N(x))\in C^2(\Omega)$ satisfies the differential inequality
\begin{equation}\label{eq0}
|Lu(x)|\leq c\sum_{|\nu|\leq 1}|\partial^{\nu}u(x)|\,\,\text{in}\,\,\Omega\,\,\text{for some}\,\,c>0,
\end{equation}
where $L$ denotes the $N\times N$ system
$$
(Lu)_{\alpha}(x)=\sum_{j,l=1}^n\sum_{\beta=1}^N C_{\alpha\beta}^{j\ell}(x)\partial_j\partial_{\ell} u_\beta(x)\,\,(1\le\alpha\le N)\,\,\text{with}\,\,
\partial_j=\partial_{x_j}\,\, (1\le j\le n)
$$
and $\partial=(\partial_1,\cdots,\partial_n)$. Through out this paper, any vector is considered as a column vector.
We first assume the following basic assumptions which are the smoothness, symmetry and strong ellipticity conditions on $L$.

\medskip
\noindent
{\sl Basic assumptions on $L$}
\begin{itemize}
\item[{\rm (i)}]
{\em smoothness condition} : Each $C_{\alpha\beta}^{j\ell}$ is a real valued Lipschitz continuous function in $\Omega$.
\item[{\rm (ii)}]
{\em symmetry condition} :
\begin{equation}\label{symm}
C_{\alpha\beta}^{j\ell}(x)=C_{\beta\alpha}^{\ell j}(x)\,\,(x\in\Omega,\,1\le j,\,\ell\le n,\,1\le\alpha,\,\beta\le N).
\end{equation}
\item[{\rm (iii)}]
{\em strong ellipticity condition} : There exists $\delta>0$ such that for any vectors
$a=(a_1,\cdots,a_N)\in\mathbb{R}^N$ and
$b=(b_1,\cdots,b_n)\in\mathbb{R}^n$
\begin{equation}\label{elliptic}
\sum_{1\le\alpha,\,\beta\le N}\sum_{1\le j,\ \ell\le n}C_{\alpha\beta}^{j\ell}(x){a}_{\alpha}b_ja_{\beta}b_{\ell}\geq\delta|a|^2|b|^2\,\,(x\in\Omega).
\end{equation}
\end{itemize}

The aim of this paper is to consider the {\em weak unique continuation
property} simply abbreviated by UCP of \eqref{eq0}, i.e., if $u(x)$ is any solution of
\eqref{eq0} which vanishes in a non-void open subset of $\Omega$,
then $u(x)$ vanishes identically in $\Omega$. We prove the UCP for
the general differential inequality (\ref{eq0}) with $L$
satisfying the above {\sl basic assumptions} and also some technical assumptions which we call {\sl further assumptions} given in the next section.
Unlike the case for scalar partial differential equations/inequalities, there are very
few results known about the UCP for systems of partial differential equations/inequalities.

A very important example of a strongly elliptic system of partial differential equations satisfying
the above three conditions is anisotropic elastic system of partial differential equations which arises
in mechanics of materials and geophysics. In this case $n=N$ and $C_{\alpha\beta}^{j\ell}$ is the elastic
tensor field. When the medium is isotropic, the UCP has been established  in \cite{aity}, \cite{dero}.
On the other hand, for the anisotropic medium, there is a rather general result known for the two dimensional case (\cite{nw}), but there are only very few results known for the three dimensional case.

A general and powerful method of proving UCP was
pioneered by Calder\'on (\cite{calderon}) and has been generalized by several authors (\cite{zuily} and
references therein).  But this is basically for scalar partial differential equations/inequalities and in order to apply this method to
systems of partial differential equations/inequalities one needs to diagonalize the leading part of system of partial differential
operators to a system of pseudo-differential operators whose characteristic roots of the principal symbol
are smooth and they are simple for real characteristic roots and at most double for complex characteristic roots.
But these conditions are usually not satisfied for general anisotropic system of elastic partial differential equations
and not even for general transversally isotropic elastic systems.

Our method, as in previous results for UCP, relies on a suitable
Carleman estimate. We adapted the idea of Vessella et al. (\cite{vess}) to reduce the derivation of a Carleman estimate to that for the operator with constant coefficients obtained by freezing the coefficients of the leading part of (\ref{eq0}) undergone the Holmgren transform. Then, we adapted the idea of Sogge (\cite{sog}) to derive a Carleman estimate for this operator with constant coefficients which is uniform with respect to the point at which we froze the coefficients. More precisely, we first
factorize this operator with constant coefficients into two first order differential operators. Next we conjugate each factor by a Carleman weight with large parameter $k$. A similar factorization was
used in \cite{na-uh} for a general elastic system of partial differential equation to study a layer stripping method for its associated inverse boundary value problem. Finally, for each conjugated factor
we construct a suitable parameterix to derive the Carleman estimate.

The rest of this paper is organized as follows. Section 2 gives some preliminaries necessary to state our main results. Especially the aforementioned {\sl further assumptions} is given there. In Section 3, we state a Carleman estimate and prove UCP using this estimate. The Carleman estimate will then be derived in the remaining sections
based on the argument outlined in the previous paragraph. More precisely, as follows. Section 4 provides the factorization of constant coefficients
system of partial differential operators $\dot{L}$ obtained by freezing the coefficients of $L$ which undergone the Holmgren transform into the product of two factors which are the first order differential operators in $t$ and first order pseudo-differential operators in the other variables. Also, they do not depend on $t$ and other variables. Here $t$ is the variable whose axis gives the direction to which we want to do UCP. Then these factors are conjugated with a Carleman weight $w_k(t)$ and obtain two pseudo-differential operators $P_{k,g}$ and $P_{k,b}$ with large parameters $k$ which do not depend on $t$ and other variables. $P_{k,g}$ is good because it is elliptic, but $P_{k,b}$ is bad because it is not elliptic. The right parametrix $S_k$ of $(T-t)^{-1}P_{k,b}^*$ is given in Section 5. The error of $S_k$ is estimated in Section 6 to show that the fundamental solution of $(T-t)^{-1}P_{k,b}^*$ can have the same estimate as that of $S_k$. The estimate of $S_k$ will be gi
 ven in Sections 7, 8. Based on these estimate for $S_k$ and that of $P_{k,g}^{-1}$, we give a Carleman estimate for $\dot{L}$ in Section 9. Finally in Section 10, we give a Carleman estimate for $L$ using the partition of unity.

\section{Preliminaries}\label{sec2}
\setcounter{equation}{0}
In order to state our main result in the next section we will give some further assumptions on the operator $L$. Let $\Gamma$ be a hypersurface in ${\Bbb R}^n$ intersecting with $\Omega$. Take $x_0\in\Gamma$ and its open neighborhood $U\subset\Omega$. Denote the symbol of $-L$ by $M=M(x,\zeta)$. For a unit conormal vector $\eta$ of $\Gamma$ and any $\xi\in{\Bbb R}^n$ independent to
$\eta$, let $\zeta=\lambda\eta+\xi$ with $\lambda\in \mathbb{R}$ and write $M$ in the form
\begin{equation}
M=T\lambda^2+A\lambda+Q,
\end{equation}
where
$$
T=T(x,\eta)=\big{(}\sum_{j,\ell=1}^n C_{\alpha\beta}^{j\ell}(x)\eta_j\,\eta_\ell;
\alpha\downarrow 1,\cdots,N,\ \beta\to 1\cdots,N\big{)},
$$
$$
A=A(x,\eta,\xi)=R(x,\eta,\xi)+R^{\top}(x,\eta,\xi),
$$
$$
R=R(x,\eta,\xi)=\big{(}\sum_{j,\ell=1}^{n} C_{\alpha\beta}^{j\ell}(x)\eta_j\xi_\ell;
\alpha\downarrow 1,\cdots,N,\ \beta\to 1\cdots,N\big{)},
$$
$$
Q=Q(x,\xi)=\big{(}\sum_{j,\ell=1}^{n} C_{\alpha\beta}^{j\ell}(x)\xi_j\xi_\ell;
\alpha\downarrow 1,\cdots,N,\ \beta\to 1\cdots,N\big{)},
$$
where the superscript "$\top$" denotes the transpose of matrices.

It is easy to see that the positivity of $T$ follows from the strong ellipticity condition of \eqref{elliptic}.
Hence $T^{1/2}$ and $T^{-1/2}$ exist and we can consider
\begin{equation}\label{factor}
H(\lambda)=H(x,\lambda,\eta,\xi)=T^{-1/2}M\,T^{-1/2}=\lambda^2+H_1(x,\eta,\xi)\lambda+H_2(x,\eta,\xi)
\end{equation}
with
\begin{equation}
H_1=T^{-1/2}AT^{-1/2}\quad\text{and}\quad
H_2=T^{-1/2}QT^{-1/2}.
\end{equation}

It is known that $H(\lambda)$ has a unique factorization
\begin{equation}\label{3.6}
H(\lambda)=(\lambda-S_1^{\ast})(\lambda-S_1)
\end{equation}
such that the spectrum $\mbox{Spec}\,S_1$ of $S_1$ satisfies Spec $S_1\subset\mathbb{C}_+:=
\{z\in\mathbb{C}\,:\,\text{\rm Im}\ z>0\}$, and this unique $S_1$ can be
given by
$$
S_1=S_1(x,\eta,\xi):=(\oint_{\gamma_+}\zeta
H(\zeta)^{-1}\,d\zeta)(\oint_{\gamma_+}H(\zeta)^{-1}d\zeta)^{-1},
$$
where $\gamma_+\subset\mathbb{C}_+:=\{\zeta\in\mathbb{C}:\ \text{\rm Im}\
\zeta>0\}$ is a closed contour enclosing all the roots of $\text{det}(H(x,\lambda,\eta,\xi))=0$ in $\lambda$ and $S_1^*$ denotes the adjoint of $S_1$ (\cite{GLR}).
Further,
\begin{equation}\label{3.7}
Z=Z(x,\eta,\xi)=\text{\rm Im}(S_1)
\end{equation}
is positive-definite (\cite{Ito}) and for $Y:=\text{Re} (S_1)$, it can be easily proved from $S_1^* S_1=H_2$, $Z^{1/2}Y  Z^{-1/2}$ is symmetric. Hence by
\begin{equation}\label{3.8}
\begin{array}{ll}
\mathcal{H}&=Z^{-1/2}H Z^{-1/2}\\
&=Z^{-1/2}\{Z^{1/2}(\lambda-Z^{-1/2}S_1^{\ast}Z^{1/2})Z^{-1/2}\}\{Z^{-1/2}
(\lambda-Z^{1/2}S_1Z^{-1/2})Z^{1/2}\}Z^{-1/2},
\end{array}
\end{equation}
$\mathcal{H}=\mathcal{H}(\lambda)$ has a much more nice factorization:
\begin{equation}\label{3.9}
\mathcal{H}(\lambda)=(\lambda-B^*)Z^{-1}(\lambda-B)
\end{equation}
with the complex symmetric matrix
\begin{equation}\label{B}
B=B(x,\eta,\xi)=B_R+iB_I, 
\end{equation}
where $B_R,\,B_I$ are given by
\begin{equation}\label{3.10}
B_R=B_R(x,\eta,\xi)=\text{Re}(Z^{1/2}S_1Z^{-1/2}),\,\,
B_I=B_I(x,\eta,\xi)=Z.
\end{equation}
Clearly the spectrums $\text{\rm Spec}(B_R\pm iB_I)$ of $B_R\pm iB_I$ satisfy
\begin{equation}\label{3.11}
\text{\rm Spec}(B_R\pm iB_I)\subset\mathbb{C}_{\pm}:=\{z\in\mathbb{C}\,:\,\pm\,\mbox{\rm Im}\,z>0\}.
\end{equation}

Now let $P_h=P_h(x,\eta,\xi)$ be the projection of $B$ given by
\begin{equation}\label{projection}
P_h=P_h(x,\eta,\xi)=(2\pi i)^{-1}\int_{\gamma_h}(\zeta-B)^{-1}d\zeta
\end{equation}
with a positively oriented small circle $\gamma_h=\gamma_h(x,\eta,\xi)$ in ${\Bbb C}_+$ enclosing the eigenvalue $\lambda_h=\lambda_h(x,\eta,\xi)\in{\Bbb C}_+$ of $B$ but excluding other eigenvalues $\lambda_k=\lambda_k(x,\eta,\xi)\in{\Bbb C}_+\,(k\not=h,\,1\le k\le s)$. Here the multiplicity $m_h$ of $\lambda_h$ satisfies $m_1+\cdots+m_s=N$. Then the further assumptions we put on $L$ are as follows.

\medskip
\noindent
{\sl Further assumptions on $L$}
\begin{itemize}
\item[{\rm (i)}] Each $P_h(x,\eta,\xi)$ is uniformly bounded for ${\rm a.e.}\,x\in U$ and all $\xi\in {\Bbb R}^n\setminus\{0\}$
linearly independent to $\eta$.
\item[{\rm (ii)}] $B(x,\eta,\xi)$ is diagonalizable for $\mbox{a.e.}\,x\in U$ and all $\xi\in {\Bbb R}^n\setminus\{0\}$
linearly independent to $\eta$.
\end{itemize}

\medskip
\begin{remark} ${}$
\newline
(a) By Theorem 1.9 in \cite{Kato}, for each almost everywhere fixed $x\in U$, the assumption {\rm (i)} implies that all the eigenvalues
$\lambda_k(x,\eta,\xi)\,(1\le k\le s)$ of $B(x,\eta,\xi)$ have to be analytic in $\xi\in{\Bbb R}^n\setminus\{0\}$ which is linearly independent to $\eta$. Needless to say this is equivalent to saying the same for all the eigenvalues of $M$.
\newline
(b) The assumption {\rm (ii)} can be given in terms of the principal symbol $M$ of $-L$. That is for $\mbox{a.e.}\,x\in U$, every $\xi\in {\Bbb R}^n\setminus\{0\}$ linearly independent to $\eta$ satisfies
\begin{equation}
\sum_{\lambda\in\Lambda(x,\xi)}\mbox{\rm dim}_{{\Bbb C}}\,\mbox{\rm Ker}\,M(x,\lambda\eta+\xi)=2N,
\end{equation}
where
\begin{equation}
\Lambda(x,\xi)=\{\lambda\in{\Bbb C}\,:\,\mbox{\rm det}\,
M(x,\lambda\eta+\xi)=0\}
\end{equation}
and $\mbox{\rm dim}_{{\Bbb C}}$ means the dimension for complex
vector spaces.
\newline
(c) Suppose we have the factorization \eqref{3.6} for \eqref{factor}. Then note that the symmetry of $Y$ is equivalent to
the commutativity of $H_1 H_2=H_2 H_1$. We refer this as {\sl commutativity}. If this commutativity holds for any $x\in U$ and $\xi\in{\Bbb R}^n\setminus\{0\}$ linearly independent to fixed $\eta$, there exist an orthogonal matrix $G=G(x,\xi,\eta)$ such that $$H_1=G\Sigma_1G^{\top},\quad H_2=G\Sigma_2G^{\top},$$
where $\Sigma_1$ and $\Sigma_2$ are diagonal matrices. Then by the strong ellipticity condition, we have $\Sigma_1^2-4\Sigma_2<0$.
Thus, we can factorize $H$ itself very nicely as follows.
$$
\begin{array}{rl}
H=&\lambda^2+\lambda G\Sigma_1G^{\top}+G\Sigma_2G^{\top}\\
=&G(\lambda^2+\lambda \Sigma_1+\Sigma_2)G^{\top}\\
=&G(\lambda+\frac{\Sigma_1}{2}-i\frac{\sqrt{4\Sigma_2-\Sigma_1^2}}{2})(\lambda+\frac{\Sigma_1}{2}+i\frac{\sqrt{4\Sigma_2-\Sigma_1^2}}{2})G^{\top}\\
=&(\lambda-B^*)(\lambda-B)\,\,\mbox{with}\,\,B=-G(\frac{\Sigma_1}{2}+i\frac{\sqrt{4\Sigma_2-\Sigma_1^2}}{2})G^{\top}.
\end{array}
$$
It is easy to check that the complex symmetric matrix $B=B_R+iB_I$ satisfies the {\sl further assumptions}. Furthermore, if we can assume $G(x,\xi,\eta)$ is smooth in $x\in U$ and $\xi\in{\Bbb R}^n\setminus\{0\}$ linearly independent to $\eta$, we can transform the operator $L$ to a system of pseudo-differential operators of order 1 which is a differential operator in $x_1$ and has diagonal principal part. By applying Calder\'on's argument(\cite{zuily}), we can have a Carleman estimate which implies the uniqueness of the Cauchy problem for the initial hyperplane $\Gamma$. If we have such a situation for every hyperplane $\Gamma$ intersecting $\Omega$, we have UCP in $\Omega$. However, this smoothness condition on $G$ does not always hold. Example 4 given below shows such a case.
\end{remark}

Next we will give several examples which satisfy all the {\sl basic assumptions} and {\sl further assumptions}. For simplicity we take $\eta=(1,0,\cdots,0)$ and write $\lambda=\xi_1,\,\xi=(\xi_2,\cdots,\xi_n)$.

\medskip
\noindent
{\bf Example 1.}
Let
\begin{align*}
Y=
\begin{pmatrix}
2 & -1\\
1 & 5
\end{pmatrix}
\xi_2,\, Z=
\begin{pmatrix}
2 & 1\\
1 & 1
\end{pmatrix}|\xi_2|.
\end{align*}
The eigenvalues of $B$ are
\begin{align*}
 \left\{
\frac{7\xi_2-\sqrt{3}|\xi_2|}{2}+\frac{3|\xi_2|+\sqrt{3}\xi_2}{2}i,\,
\frac{7\xi_2+\sqrt{3}|\xi_2|}{2}+\frac{3|\xi_2|-\sqrt{3}\xi_2}{2}i
\right\}
\end{align*}
and the corresponding eigenvectors are
\begin{align*}
& (10|\xi_2|+2\xi_2i,\,
 5\sqrt{3}\xi_2-5|\xi_2|+(5\sqrt{3}|\xi_2|-11\xi_2)i),\\
& (-10|\xi_2|-2\xi_2i,\, 5\sqrt{3}\xi_2+5|\xi_2|+(5\sqrt{3}|\xi_2|+11\xi_2)i).
\end{align*}
The associated $H$ is
\begin{align*}
H=\xi_1^2I+
\begin{pmatrix}
-4 & 0\\
0 & -10
\end{pmatrix}
\xi_1\xi_2+
\begin{pmatrix}
10 & 6\\
 6 & 28
\end{pmatrix}
\xi_2^2.
\end{align*}
and it does not have the commutativity.

\medskip
\noindent
{\bf Example 2.} Perturb the previous example as follows.
Let  \begin{align*}
Z=
\begin{pmatrix}
2 & 1\\
1 & 1
\end{pmatrix}
|\xi_2|+
\begin{pmatrix}
q\varepsilon_1 & 0\\
0 & q\varepsilon_2
\end{pmatrix}
|\xi_2|
\end{align*}
with $\varepsilon_1$, $\varepsilon_2$, $q\in \mathbb{R}$ satisfying
$|q+\frac{3}{\varepsilon_1+\varepsilon_2}|\geq\delta>0$
and $|q-\frac{1}{\varepsilon_1-\varepsilon_2}|\geq\delta>0$ for some
$\delta$.
The eigenvalues of $B$ are
\begin{align*}
&\lambda_{\pm}\\
:=&\frac{14\xi_2\mp c|\xi_2|}{4}+\frac{i}{2(1+a)}
\left((3+4q\varepsilon_1+q^2\varepsilon^2_1-2q\varepsilon_2-q^2\varepsilon_2^2)
|\xi_2|\pm\frac{c\xi_2}{12}(\sqrt{b}+a^2+2a)\right),
\end{align*}
where
\begin{align*}
 a=q(\varepsilon_1-\varepsilon_2),\, b=a^4+4a^3+40a^2+72a+36,\,
 c= \sqrt{2}\sqrt{\sqrt{b}-a^2-2a}.
\end{align*}
Corresponding eigenvectors are
\begin{align*}
 (\beta,\, \lambda_{\pm}-\alpha),
\end{align*}
where
\begin{align*}
& \beta =
 \frac{(d^2+d(2q\varepsilon_1+1)+q\varepsilon_1-1)\xi_2+i|\xi_2|}
{d(3+q(\varepsilon_1+\varepsilon_2)+2d^2)},\\
& \alpha
 = \frac{(4d^2+q(\varepsilon_1+\varepsilon_2))\xi_2+i(2+q\varepsilon_1)|\xi_2|}
{d(3+q(\varepsilon_1+\varepsilon_2)+2d^2)},\\
 & d=\sqrt{1+q\varepsilon_1+2q\varepsilon_2+q^2\varepsilon_1\varepsilon_2}.
\end{align*}
The associated $H$ is
\begin{align*}
H=\xi_1^2I+
\begin{pmatrix}
-4 &0 \\
0 & -10
\end{pmatrix}
\xi_1\xi_2+
\begin{pmatrix}
10+4q\varepsilon_1+q^2\varepsilon_1^2 &  6+q\varepsilon_1+q\varepsilon_2\\
6+q\varepsilon_1+q\varepsilon_2       &  28+2q\varepsilon_2+q^2\varepsilon_2^2
\end{pmatrix}
\xi_2^2.
\end{align*}
This $H$ also does not have the commutativity.

\medskip
\noindent
{\bf Example 3.}
Let
\begin{align*}
& Y=
\begin{pmatrix}
 1 & 0 & 0\\
 0 & 3 & 1 \\
 0 & -1& 0
\end{pmatrix}
\xi_2+
\begin{pmatrix}
 1 & 0 & 0\\
 0 & -1 & 1 \\
 0 & -1 & -4
\end{pmatrix}
\xi_3,\\
& Z=
\begin{pmatrix}
 1 & 0 & 0\\
 0 & 2 & 1 \\
 0 & 1 & 1
\end{pmatrix}
\sqrt{3\xi_2^2+2\xi_2\xi_3+\xi_3^2}.
\end{align*}
The eigenvalues of $B$ are
\begin{align*}
& \lambda_1:=\xi_2+\xi_3+i\sqrt{3\xi_2^2+2\xi_2\xi_3+\xi_3^2},\\
& \lambda_{\pm}:=\frac{3\xi_2-5\xi_3\pm b}{2}+
i\frac{9(3\xi_2^3+5\xi_2\xi_3^2+3\xi_2^2\xi_3+\xi_3^3)\pm b\sqrt{a}\pm
 5\xi_2^2b}{6(\xi_2+\xi_3)\sqrt{3\xi_2^2+2\xi_2\xi_3+\xi_3^2}},
\end{align*}
where
\begin{align*}
 b:= \sqrt{\sqrt{a}-5\xi_2^2},\,
 a:= 52\xi_2^4+72\xi_2^3\xi_3+72\xi_2\xi_3^2+36\xi_2\xi_3^3+9\xi_3^4.
\end{align*}
The corresponding eigenvectors are
\begin{align*}
 (1,0,0),\, (\beta,\, \lambda_{\pm}-\alpha),
\end{align*}
where
\begin{align*}
 \beta:=\frac{\xi_2+\xi_3}{5}+i\sqrt{3\xi_2^2+2\xi_2\xi_3+\xi_3^2},\,
\alpha:= \frac{13\xi_2-7\xi_3}{5}+2i\sqrt{3\xi_2^2+2\xi_2\xi_3+\xi_3^2}.
\end{align*}
The associated $H$ is
\begin{align*}
H=&\xi_1^2I+\xi_1
\begin{pmatrix}
-2\xi_2-2\xi_3& 0& 0\\
0 &-6\xi_2+2\xi_3& 0 \\
0 & 0 & 8\xi_3
\end{pmatrix}
\\
&+
\begin{pmatrix}
4\xi_2^2+4\xi_2\xi_3+6\xi_3^2& 0 & 0 \\
0 &  22\xi_2^2+4\xi_2\xi_3+6\xi_3^2  &  6\xi_2^2+8\xi_2\xi_3+4\xi_3^2 \\
0 &  6\xi_2^2+8\xi_2\xi_3+4\xi_3^2  &  4\xi_2^2+4\xi_2\xi_3+18\xi_3^2
\end{pmatrix}.
\end{align*}
This $H$ does not have the commutativity.

%\begin{remark}
%The eigenvalues and the coresponding eigenvectors of $B$ 
%in Example 2 and Example 3 are complicated
%expressions. However, we can easily obtain them from ${\rm Re}\, S_1$ and
%$Z$. The eigenvalues are the solutions of quadratic equations with
%complex coefficients, therefore we can obtain the explicit formula of
% the eigenvalues and the corresponding eigenvectors.
%\end{remark}
%\medskip
%\noindent
{\bf Example 4.} Let
\begin{equation}
H=(\xi_1-B^*)(\xi_1-B)\,\,\text{with}\,\,B=B_R+i B_I,
\end{equation}
where
\begin{equation*}
B_R=
\begin{pmatrix}
0         &     0      & x_2\xi_2        \\
0         &     0      & x_3\xi_3        \\
x_2\xi_2  &  x_3\xi_3   & 0
\end{pmatrix}, \quad
B_I=\sqrt{\xi_2^2+\xi_3^2}\,I
\end{equation*}
The eigenvalues of $B$ are
$$
\left\{
 -\sqrt{(x_2\xi_2)^2 + (x_3\xi_3)^2} + i\sqrt{\xi_2^2+\xi_3^2},\,
 \sqrt{(x_2\xi_2)^2 + (x_3\xi_3)^2} + i\sqrt{\xi_2^2+\xi_3^2},\,
 i\sqrt{\xi_2^2+\xi_3^2}
 \right\},
$$
and the corresponding eigenvectors outside the set $\{x_2=x_3 = 0\}$ are given by
$$
v_1 := \left(x_2\xi_2,\, x_3\xi_3,\, - \sqrt{(x_2\xi_2)^2 + (x_3\xi_3)^2}\right),
$$
$$
v_2 := \left(x_2\xi_2,\, x_3\xi_3,\, \sqrt{(x_2\xi_2)^2 + (x_3\xi_3)^2}\right)
$$
and
$$
v_3 := \left(x_3\xi_3,\, -x_2\xi_2,\, 0\right),
$$
respectively.
This example shows that all the conditions are satisfied if $x_2 x_3\not=0$
and it does have the commutativity.
However, there is no continuous diagonalization where we have
$x_2=x_3=0$. As a matter of fact, the $1$-dimensional vector bundle
spanned by the vector $v_3$ on $\{x^2_2 + x^2_3 \ne 0\}$ never extends
to the whole space continuously since
the point $(x_3\xi_3,\,-x_2\xi_2)$
turns once around the origin in $\mathbb{R}^2$ if $\xi_2\xi_3 \ne 0$
and $(x_2,\,x_3)$ moves once around the origin.

\medskip
\noindent
{\bf Example 5.} Let
\begin{equation*}
\begin{array}{rl}
H=&\big[\,\xi_1I-(A^{-1/2}EA^{1/2})^{\top}+iA\sqrt{\xi_2^2+\xi_3^2}\,\big]\\
&\times
\big[\,\xi_1I-A^{-1/2}EA^{1/2}-iA\sqrt{\xi_2^2+\xi_3^2}\,\big]\\
=&\xi_1^2I-(A^{-1/2}EA^{1/2})^{\top}-(A^{-1/2}EA^{1/2})\\
&+(A^{-1/2}EA^{1/2})^{\top}(A^{-1/2}EA^{1/2})+(\xi_2^2+\xi_3^2)A^2,
\end{array}
\end{equation*}
where $A=
\begin{pmatrix}
1  &  0  &  0  \\
0  &  1+ (x_2x_3)^2  & 0        \\
0  &  0 &  1- (x_2x_3)^2
\end{pmatrix}$ and $E=
\begin{pmatrix}
1     &      0     &  0   \\
0     &  x_2\xi_2  & x_3\xi_3        \\
0     &   x_3\xi_3 &  -x_2\xi_2
\end{pmatrix}$.
We assume that, in what follows, $(x_2,\,x_3)$ is sufficiently close to the origin.
The eigenvalues of $S_1=Y+iZ=A^{-1/2}EA^{1/2}+iA\sqrt{\xi_2^2+\xi_3^2}$ which is in fact equal to $B=B_I+iB_I$ are given by $ir_\xi$,
$$
ir_\xi-\sqrt{ - r_\xi^2 x_2^4 x_3^4+2\,i\,\xi_2\,
r_\xi\, x_2^{3}\,x_3^{2}+r_{x\xi}^2}
$$
and
$$
ir_\xi+\sqrt{ - r_\xi^2 x_2^4 x_3^4+2\,i\,\xi_2\,
r_\xi\, x_2^{3}\,x_3^{2}+r_{x\xi}^2},
$$
where we set $r_\xi := \sqrt{\xi_2^2 + \xi_3^2}$ and
$r_{x\xi} := \sqrt{(x_2\xi_2)^2 + (x_3\xi_3)^2}$.
The corresponding eigenvectors on $x_2x_3\xi_2\xi_3 \ne 0$ are
$v_1 :=(1,\,0,\,0)$,
$$
v_2 := (0,\,\,
-\sqrt{-r_{\xi}^2  x_2^4 x_3^4
+ 2\,i\,\xi_2\,r_\xi x_2^3 x_3^2 + r_{x\xi}^2}+
i\,r_{\xi} x_2^{2} x_3^{2}+ \xi_2 x_2,\,\,
\xi_3 x_3)
$$
and
$$
v_3 := (0,\,\,
\sqrt{-r_{\xi}^2  x_2^4 x_3^4
+ 2\,i\,\xi_2\,r_\xi x_2^3 x_3^2 + r_{x\xi}^2}+
i\,r_{\xi} x_2^{2} x_3^{2}+ \xi_2 x_2,\,\,
\xi_3 x_3),
$$
repsectively. 
%Note also that $v_2$ is parallel to the vector
%$$
%(0,\,\,
%-\xi_3 x_3,\,\,
%\sqrt{-r_{\xi}^2  x_2^4 x_3^4
%+ 2\,i\,\xi_2\,r_\xi x_2^3 x_3^2 + r_{x\xi}^2}+
%i\,r_{\xi} x_2^{2} x_3^{2}+ \xi_2 x_2),
%$$
%and $v_3$ is parallel to the vector
%$$
%(0,\,\,
%-\xi_3 x_3,\,\,
%-\sqrt{-r_{\xi}^2  x_2^4 x_3^4
%+ 2\,i\,\xi_2\,r_\xi x_2^3 x_3^2 + r_{x\xi}^2}+
%i\,r_{\xi} x_2^{2} x_3^{2}+ \xi_2 x_2).
%$$
This example shows that all the conditions are satisfied
if $x_2 x_3 \not=0$ and it does not have the commutativity.
Also, there is no continuous diagonalization where we have $x_2=x_3=0$.
These observations follow from the fact that there exists a constant $M > 0$ such that,
on $\xi_2^2 + \xi_3^2 = 1$, we have
$$
\left|v_2 - \left(0,\,\, - r_{x\xi}+ \xi_2 x_2,\,\, \xi_3 x_3\right)\right|
\le M r_{x\xi}^2
$$
and
$$
\left|v_3 - \left(0,\,\, r_{x\xi} + \xi_2 x_2,\,\, \xi_3 x_3\right)\right|
\le M r_{x\xi}^2.
$$

\vskip 1cm
Through out this paper we will always assume the above {\sl basic assumptions} and {\sl further assumptions}.

\section{Main results}

First of all as it has been already mentioned before in the last remark given in the previous section, the UCP in $\Omega$ follows from the following uniqueness of the Cauchy problem (see \cite{zuily}) for any hyperplane intersecting with $\Omega$. Also, the uniqueness of the Cauchy problem follows from a suitable Carleman estimate given below in this section.
\begin{theorem}\label{uni-cau}
Let $\Gamma$, $x_0$, $U$ be as in the
{\sl further assumptions}. Also let $u\in C^2(U)$ satisfy
\begin{equation}\label{cauchyp}
\begin{cases}
|Lu|\leq c\sum_{|\nu|\leq 1}|\partial^{\nu}u|\quad\text{in}\quad U\\
{\partial}^{\mu}u|_{\Gamma}=0\,\, (|\mu|\leq 1).
\end{cases}
\end{equation}
Then there exists a small open neighborhood of $U'\subset U$ of $x_0$ so that
$u\equiv 0$ in $U'$.
\end{theorem}

Next we will describe the Carleman estimate which we will use to prove Theorem \ref{cauchyp}.
Let $\eta^0$ be the unit co-normal vector of $\Gamma$. By an appropriate translation
and rotation, we can assume that $x_0=0$ and
$\eta^0=(1,0,\cdots,0)^{\top}$. We now make a change of
coordinates near $0$ by using the ``Holmgren transform", i.e.,
$$x_j\to x_j,\quad 2\leq j\leq n,\quad t=x_1+\kappa|x'|^2\quad\text{for an appropriate constant}\ \kappa,$$ where
$x'=(x_2,\cdots,x_n)$. In the new coordinates
$\widetilde{x}=(\tilde{x}^1,\tilde{x}')=(\tilde{x}_1,\tilde{x}_2,\cdots,\tilde{x}_n)=(t,x')$, \eqref{cauchyp} becomes
\begin{equation}\label{cauchyn}
\begin{cases}
|\widetilde{L}\widetilde u|\leq\widetilde{c}\sum_{|\nu|\leq 1}|D_{\widetilde{x}}\widetilde{u}|\quad\text{in}\quad \widetilde{U}\\
{\rm supp}(\widetilde{u})\subset\{\widetilde{x}:
t\geq\widetilde{\kappa}|\tilde{x}'|^2\}\quad\text{for some constants}\
\widetilde c,\widetilde{\kappa}>0,
\end{cases}
\end{equation}
where $\widetilde{u}(\widetilde{x})=u(x(\widetilde{x}))$,
$\widetilde{U}$ is a small neighborhood of $0$ and the partial
differential operator $\widetilde{L}$ is defined by
\begin{equation}\label{tildeL}
(\widetilde{L}\widetilde{u})_{\alpha}=\sum_{1\le p,\,q\le n}\sum_{1\le\beta\le N}\widetilde{C}_{\alpha\beta}^{pq}(\widetilde{x})\partial_p\partial_q\widetilde{u}_{\beta}
\end{equation}
with
$$
\widetilde{C}_{\alpha\beta}^{pq}(\widetilde{x})=J^{-1}(x(\widetilde{x}))\sum_{1\le j,\,\ell\le n}C_{\alpha\beta}^{j\ell}(x(\widetilde{x}))\partial_j\widetilde{x}_p\partial_{\ell}\widetilde{x}_q,
$$
the Jacobian $J(x(\widetilde{x}))$ of the change of variables and partial derivative $\partial_p$ with respect to $\tilde{x}_p$.
It can be easily checked that $\widetilde{L}$ satisfies
the {\sl basic assumptions} and {\sl further assumptions} in $\widetilde U$.
\begin{theorem}\label{thm2.2}
There exist positive constants $T_0$, $k_0$, $r$, and $c$ such
that for $0<T\leq T_0$, we have for $k\geq k_0$ that
\begin{equation}\label{2.3}
\sum_{|\nu|\leq 1}T^{-1}(kT)^{3-2|\nu|}\int_0^Te^{k(t-T)^2}\|\partial^{\nu}v\|^2dt\leq
c\int_0^Te^{k(t-T)^2}\|\widetilde{L}v\|^2dt
\end{equation}
holds for all $v(t,\tilde{x}')\in C^{\infty}({\Bbb R}^n)$ with $\ {\rm supp}(v)\subset\{(t,
\tilde{x}'):
|\tilde{x}'|< r,\ 0< t< T/2\}$, where $\|\cdot\|^2=(\cdot,\cdot)$
is the $L^2(\mathbb{R}^{n-1})$ norm.
\end{theorem}

Theorem~\ref{thm2.2} will be proved in the remaining sections. Once having this
Carleman estimate,
the proof of Theorem~\ref{uni-cau}
is rather standard (see for instance
\cite{zuily}). We give the proof here for the sake of
completeness.

\medskip{\em Proof of Theorem~\ref{uni-cau}}. Let
$\theta(t)\in C^{\infty}$ defined in $t\geq 0$ with $\theta(t)=0$
for $t\geq T$ and $\theta(t)=1$ for $t\leq 2T/3$. Assume that
$\widetilde{u}$ is a solution of \eqref{cauchyn}. When $T$ is
sufficiently small, we can apply \eqref{2.3} to
$\theta\widetilde{u}$ and get that
\begin{equation}\label{2.4}
\begin{array}{rl}
\sum_{|\nu|\leq
1}T^{-1}(kT)^{3-2|\nu|}\int_0^Te^{k(t-T)^2}\|\partial^{\nu}(\theta\widetilde{u})\|^2dt&\leq
c\int_0^Te^{k(t-T)^2}\|\widetilde{L}(\theta\widetilde{u})\|^2dt\\
{}&\leq
c\int_0^Te^{k(t-T)^2}\|\theta\widetilde{L}\widetilde{u} \|^2dt\\
{}&\quad+c\int_0^Te^{k(t-T)^2}\|\
[\widetilde{L},\theta]\widetilde{u}\|^2dt,
\end{array}
\end{equation}
where $[\widetilde{L}, \theta]$ denotes the commutator of $\widetilde{L}$ and $\theta$.
By using the equation in \eqref{cauchyn} and taking $k$ large, $T$
small, if necessary, we can absorb the first term on the right
side of \eqref{2.4}. Hence we have
\begin{equation}\label{2.5}
\sum_{|\nu|\leq
1}\int_0^{T/2}e^{k(t-T)^2}\|\partial^{\nu}\widetilde{u}\|^2dt\leq
c'\sum_{|\nu|\leq
1}\int_{2T/3}^Te^{k(t-T)^2}\|\partial^{\nu}\widetilde{u}\|^2dt
\end{equation}
for some $c'>0$. Since $e^{k(t-T)^2}$ is decreasing in $0\leq
t\leq T$, we have from \eqref{2.5} that
$$
\sum_{|\nu|\leq
1}\int_0^{T/2}\|\partial^{\nu}\widetilde{u}\|^2dt\leq
c'e^{k(-T^2/4+T^2/9)}\sum_{|\nu|\leq
1}\int_{2T/3}^T\|\partial^{\nu}\widetilde{u}\|^2dt,
$$
which implies that $\widetilde{u}\equiv 0$ in $t\leq T/2$. \eproof

\section{Associated constant coefficients operators}
\setcounter{equation}{0}
In this section we provide some factorization of $\widetilde{L}$
and consider an associated constant coefficient operator. To avoid any further heavy notations, we suppress using "$\,\,\widetilde{}$\,\,"\,in this sections but also in the rest of sections except the last section. That is we abuse the notations $L$ etc. to denote $\widetilde{L}$ etc.

From the {\sl basic assumptions}, we have the factorization (\ref{3.9}) with simplified expression on the dependency of $\eta,\,\xi$. That is by taking $\eta=(1,0,\cdots,0)$, $\xi=(\xi_1,\xi')=(\lambda,\xi')$ with $\xi'=(\xi_2,\cdots,\xi_n)\in{\Bbb R}^{n-1}$, we can
simplify the notations, for example $M$ and $P_h$ can be written as $M=M(x,\xi)$ and $P_h=P_h(x,\xi')$, respectively. Since it is enough to obtain a Carleman estimate for $T^{1/2} L T^{1/2}$, we simply denote it by $L$.

On this occasion, we introduce the notations $D=-i(\partial/(\partial x_1),\cdots,\partial/(\partial x_n))$ and $D'=-i(\partial/(\partial x_2),\cdots,\partial/(\partial x_n))$.
Further we note here that we sometimes write $t$ instead of $x_1$, that is $t=x_1$.

Now we conjugate $\mathcal{L}:=-L$ by the Carleman weight function $w_k(t):=\mbox{exp}(k(t-T)^2/2)$ with a large parameter $k>0$. That is to consider
\begin{equation}\label{tildeLx}
\mathcal{L}_k:=w_k(t)\mathcal{L}w_k(t)^{-1}.
\end{equation}
The principal symbol of $\mathcal{L}_k$ under the scaling
\begin{equation}\label{scaling}
(\lambda,\xi',k)\rightarrow(\rho\lambda,\rho\xi',\rho k)\,\,
\mbox{\rm for}\,\,\rho>0
\end{equation}
is given by
\begin{equation}\label{PrincipalSymbol}
H_k=H_k(x,\xi,k)=\{\lambda-[Y^{\top}-i(Z+k(t-T))]\}\{\lambda-[Y+i(Z-k(t-T))]\}.
\end{equation}

We first fix $x_0\in U$ and try to derive a Carleman estimate
for $\dot{\mathcal{L}}$ obtained from $\mathcal{L}$ by freezing its coefficients at $x_0$. Note that $\dot{\mathcal{L}}$ is given by
\begin{equation}
\begin{array}{l}
\dot{\mathcal{L}}(D_t,D')
=\{D_t-(\dot{Y}(D')^{\top}-i\dot{Z}(D'))\}\{D_t-(\dot{Y}(D')+\dot{Z}(D'))\},
\end{array}
\end{equation}
where $\dot{Y}(\xi')=Y(x_0,\xi')$, $\dot{Z}(\xi')=Z(x_0,\xi')$.
By defining $\dot{H}(\xi)=H(x_0,\xi)$, we have
\begin{equation}
A_k=A_k(D_t,D')=w_k(t)\dot{H}(D_t,D')w_k(t)^{-1}=P_{k,b}P_{k,g},
\end{equation}
where the good part $P_{k,g}$ and bad part $P_{k,b}$ are given by
\begin{equation}
P_{k,g}=D_t-\dot{G}_k(D'),\,\,P_{k,b}=D_t-\dot{G}_k^*(D')
\end{equation}
with
\begin{equation}
\dot{G}_k(\xi')=\dot{Y}(\xi')+i\big(\dot{Z}(\xi')-k(t-T)\big),\,\,\,
\dot{G}_k^*(\xi')=\dot{Y}(\xi')^{\top}-i\big(\dot{Z}(\xi')+k(t-T)\big).
\end{equation}

   To have an a priori estimate for $P_{k,g}$ and $P_{k,b}$, we construct right parametrices $\mathcal{S}_{k,g}$ and $\mathcal{S}_{k,b}$ for their
adjoint operators $P_{k,g}^*$ and $P_{k,b}^*$, respectively. Since $P_{k,g}^*$ is an elliptic operator of order one with large parameter $k$,
$\mathcal{S}_{k,g}$ can be easily constructed using the theory of pseudo-differential operator with large parameter to have
\begin{equation}\label{remainderR(k,g)}
P_{k,g}^*\mathcal{S}_{k,g}=I+R_{k,g}
\end{equation}
for some pseudo-differential operator $R_{k,g}$ of large negative order with large parameter $k$ which satisfies the estimate
\begin{equation}\label{estimageR(k,g)}
\Vert (I+R_{k,g})^{-1}\Vert\le1\,\,\,(kT\ge M_g)
\end{equation}
for some $M_g>0$, where the norm $\Vert\cdot\Vert$ is the operator norm on $L^2({\Bbb R}^n)$.
Further, for each $s\in{\Bbb R}$, there exists a constant $C_g(s)>0$ such that $\mathcal{S}_{k,g}$ satisfies the estimate
\begin{equation}\label{estimateS(k,g)}
\Vert \mathcal{S}_{k,g} w\Vert_{H^s({\Bbb R}^n)}\le C_g(s)(kT)^{-1}\Vert w\Vert_{H^s({\Bbb R}^n)}\,\,\,(w\in H^s({\Bbb R}^n)).
\end{equation}

As for $\mathcal{S}_{k,b}$ we construct a right pararametrix $\mathcal{S}_k':=\mathcal{S}_{k,b}'$ for
\begin{equation}
\mathcal{P}_k^*=\mathcal{P}_{k,b}^*=\frac1{T-t}P_{k,b}^*
\end{equation}
in the next section.

\section{A right parametrix $\mathcal{S}_k'$ of $\mathcal{P}^*_k$}\label{sec5}
We look for the Schwartz kernel $S_k'(x_1,y_1;\xi')$ of a right parametrix $\mathcal{S}_k'$ of the ordinary differential operator $\mathcal{P}_k^*$ in $x_1$ with parameters $\xi',\,k$ near its characteristics in the form
\begin{equation}\label{5.2}
S_k'(x_1,y_1;\xi')=\int \eta(\frac{\xi_1}{2\gamma kT}) e^{i\Phi(x_1,y_1,\xi')}M^{-1}(y_1,\xi')
\tilde{d}\xi_1,
\end{equation}
where
\begin{equation}\label{5.3}
\begin{array}{rl}
\Phi(x_1,y_1,\xi')&=(x_1-y_1)\xi_1+\Phi_1(x_1,y_1,\xi')\\
M(y_1,\xi')&=\frac{1}{T-y_1}\Big(\xi_1-\dot{Y}(\xi')-i\big(\dot{Z}(\xi')-k(T-y_1)\big)\Big)
\end{array}
\end{equation}
and $\tilde{d}\xi_1=(2\pi)^{-1}d\xi_1$,  $\eta\in C^\infty_0(\mathbb{R})$, $\eta(s)=1$ for $|s|\leq 1/4$ and $\eta(s)=0$ for $|s|\geq 1/2$ with $0\leq\eta\leq1$.
Further $\gamma>0$ is defined as follows.
The principal symbol of $\sigma(P_k^*)$ of $P_k^*$ is given by
\begin{equation}\label{5.4}
\begin{array}{ll}
\sigma(P_k^*)&=\xi_1-\dot{Y}(\xi')-i\big(\dot{Z}(\xi')-k(T-x_1))\big)\\
&=\dot{Z}^{-1/2}(\xi')\{\xi_1-\big(\dot{B}(\xi')-ik(T-x_1))\big)\}\dot{Z}^{1/2}(\xi'),
\end{array}
\end{equation}
where $\dot{B}=B(x_0,\xi')$ is the simplified notation of $B(x_0,\eta,\xi)$ given by \eqref{B} for $\eta=(1,0,\cdots,0)$, $\xi=(\xi_1,\xi')=(\lambda,\xi')$.

Due to the symmetry of  $B_R(x,\xi')$ and $B_I(x,\xi')$, observe that for any $v\in\mbox{\rm Ker}\,\sigma(P_k^*)$ with $|v|=1$, we have
\begin{equation}\label{5.5}
((\dot{B}_I(\xi')-k(T-x_1))v,v)=0
\end{equation}
and
\begin{equation}\label{5.6}
((\xi_1-\dot{B}_R(\xi'))v,v)=0.
\end{equation}
Further since $B_I(x,\xi')=Z(x,\xi')$ is positive, we have for $x_1\in[0,T/2]$,
\begin{equation}\label{5.7}
\lambda^{-1}(kT)\le|\xi'|\le\lambda(kT)
\end{equation}
for some constant $\lambda>0$ from \eqref{5.5} and
\begin{equation}\label{5.8}
|\xi_1|\le m|\xi'|\,\,\mbox{\rm with}\,\,m=\mbox{\rm max}\{|\mbox{\rm inf}_{|v|=1}(\dot{B}_Rv,v)|,\,
|\mbox{\rm sup}_{|v|=1}(\dot{B}_Rv,v)|\}
\end{equation}
from \eqref{5.6}. Combining \eqref{5.8} with \eqref{5.7},
we have
\begin{equation}\label{5.9}
|\xi_1|\le \gamma(kT)\,\,\text{with}\,\,\gamma=\lambda m.
\end{equation}
Note that $m$ can be taken uniformly for $x_0\in U$.

Based on these the characteristics of $\mathcal{P}_k^*$ is simply given by
\begin{equation}
\text{det}\{\xi_1-\big(\dot{B}(\xi')-ik(T-x_1)\big)\}=0
\end{equation}
for large $k>0$ and it is enough to consider the operator with prinicipal symbol
$\xi_1-\big(\dot{B}(\xi')-ik(T-x_1)\big)$ and $(T-y_1)^{-1}\big(\xi_1-\big(\dot{B}(\xi')-ik(T-y_1)\big)$ instead of $\mathcal{P}_k^*$
and $M(y_1,\xi')$, respectively. In the rest of this paper, we denote them by $\mathcal{P}_k^*$ and $M(y_1,\xi')$, respectively.

Since we want to have
\begin{equation}\label{5.10}
\begin{array}{rl}
I&=e^{-i\Phi(x_1,y_1,\xi')}\mathcal{P}_k^*(e^{i\Phi(x_1,y_1,\xi')})M^{-1}(y_1,\xi')\\
&=[\frac1{T-x_1}\Phi_{1,x_1}+\frac{\xi_1}{T-x_1} -\frac1{T-x_1}\dot{B}_R(\xi')-i\frac1{T-x_1}\dot{B}_I(\xi')+ik] M^{-1}(y_1,\xi'),
\end{array}
\end{equation}
we can take
\begin{equation}\label{5.11}
\Phi_1(x_1,y_1,\xi')=\frac{(x_1-y_1)^2}{2(T-y_1)}(-\xi_1+\dot{B}_R+i\dot{B}_I).
\end{equation}

Now let
\begin{equation}\label{5.12}
\begin{array}{ll}
S_{k,1}'(x,y)&=\int_{{\Bbb R}^{n-1}}\theta_\lambda(\frac{\xi'}{kT})e^{i(x'-y')\cdot\xi'}\Big(\int \eta(\frac{\xi_1}{2\gamma kT}) e^{i\Phi(x_1,y_1,\xi')}M^{-1}(y_1,\xi')\tilde{d}\xi_1\Big)\tilde{d}\xi'\\
&=\int_{{\Bbb R}^{n-1}}\theta_\lambda(\frac{\xi'}{kT})e^{i(x'-y')\cdot\xi'}S_k'(x_1,y_1;\xi')\tilde{d}\xi'
\end{array}
\end{equation}
and
\begin{equation}\label{5.13}
\begin{array}{rl}
S_{k,0}'(x,y)&=\int_{{\Bbb R}^{n-1}}e^{i(x'-y')\cdot\xi'}\Big(\int \beta_0(\xi) e^{i(x_1-y_1)\xi_1}M^{-1}(y_1,\xi')\tilde{d}\xi_1\Big)\tilde{d}\xi'\\
&=\int_{{\Bbb R}^{n}}e^{i(x-y)\cdot\xi}\beta_0(\xi) M^{-1}(y_1,\xi')\tilde{d}\xi,
\end{array}
\end{equation}
where
$\tilde{d}\xi'=(2\pi)^{-(n-1)}d\xi'$, $\tilde{d}\xi=(2\pi)^{-n)}d\xi$ and $\theta_\lambda(\xi')\in C_0^\infty({\Bbb R}^{n-1})$ such that
$\theta_\lambda(\xi')=1\,((2\lambda)^{-1}\le|\xi'|\le 2\lambda)$, $\theta_\lambda(\xi')=0\,(\xi'\sim0)$ and
$\beta_0=1-\beta_1$ with $\beta_1(\xi)=\theta_\lambda(\xi'/(kT))\eta(\xi_1/(2\gamma kT))$.
Then, by decomposing $S_k'(x,y)$ into $S_{k}'(x,y)=S_{k,1}'(x,y)+S_{k,0}'(x,y)$, we will prove that $S_k'(x,y)$ is the Schwartz kernel of a right parametrix $\mathcal{S}_k'$ of $\mathcal{P}_k^*$ in the next section.

Before closing this section, we remark here that we have
\begin{equation}\label{5.14}
|\mbox{\rm det}\Big(\beta_0(\xi) M(y_1,\xi')\Big)|\gtrsim T^{-N}(|\xi|+kT)^N,
\end{equation}
where the notation "$\gtrsim$" denotes "$\ge$" modulo  multiplication by a general positive constant.

\section{Error estimate of parametrix $\mathcal{S}_k'$}\label{sec6}
In order to see that $\mathcal{S}_k'$ is a parametrix of $\mathcal{P}_k^*$, we need to estimate the error
\begin{equation}\label{7.1}
\mathcal{P}_k^* \mathcal{S}_k'-I.
\end{equation}
Since $\int\beta_0(\xi)e^{i(x-y)\cdot\xi}\tilde{d}\xi+\int\beta_1(\xi)e^{i(x-y)\cdot\xi}\tilde{d}\xi=\delta(x-y)$ with $\tilde{d}\xi=(2\pi)^{-n}d\xi$ and taking into account (\ref{5.10}), the sum of two errors coming from the operators with Schwartz's kernels $S_{k,1}'(x,y),\, S_{k,0}'(x,y)$ is
\begin{equation}\label{7.2}
\begin{array}{l}
\int_{{\Bbb R}^{n-1}}e^{i(x'-y')\cdot\xi'}\{\int \beta_1(\xi)
 e^{i\Phi_I(x_1,y_1,\xi')}\tilde{d}\xi_1\}\tilde{d}\xi'\\
 +\int_{{\Bbb R}^{n-1}}e^{i(x'-y')\cdot\xi'}\{\int \beta_0(\xi)
 e^{i(x_1-y_1)\xi_1}\tilde{d}\xi_1\}\tilde{d}\xi'-\delta(x-y)
\end{array}
\end{equation}
in terms of Schwartz kernel.

Next we will show that by taking $k$ and $T$ to satisfy
$k\ge T^{-3}$, we have
\begin{equation}\label{7.3}
\begin{array}{l}
\int\beta_0(\xi) e^{i(x_1-y_1)\xi_1}\tilde{d}\xi_1 + \int \beta_1(\xi)e^{i\Phi(x_1,y_1,\xi')} \tilde{d}\xi_1\\
=\delta(x_1-y_1)+O(k^{-1/6})
\end{array}
\end{equation}
over a neighborhood of $x_1=y_1$ uniformly respect to any fixed $\xi'\in{\Bbb R}^{n-1}$, $y_1\in[0,T/2]$, where $O(k^{-1/6})$ is the kernel of a bounded operator on $L^2((0,T/2))$
norm with operator estimated by $O((k)^{-1/6})$ uniformly with respect to other variables such as $y_1, \xi'$.
For this it is enough to show
\begin{equation}\label{7.4}
\begin{array}{ll}
\int \beta_1(\xi) e^{i\Phi(x_1,y_1,\xi')} \tilde{d}\xi_1
&=\int e^{i\Phi_2(x_1,y_1,\xi_1)}\beta_1(\xi) e^{\Phi_3(x_1,y_1,\xi')} \tilde{d}\xi_1\\
&=\int\beta_1(\xi) e^{i(x_1-y_1,\xi_1)} \tilde{d}\xi_1 +O((kT^2)^{-1})
\end{array}
\end{equation}
over a neighborhood of $x_1=y_1$, where $\Phi_2(x_1,y_1,\xi_1)=(x_1-y_1)\xi_1-\frac{(x_1-y_1)^2}{2(T-y_1)}\xi_1$
and $\Phi_3(x_1,y_1,\xi')=\frac{i(x_1-y_1)^2}{2(T-y_1)}\dot{B}(\xi')$. Write
\begin{equation}\label{7.5}
\begin{array}{l}
\Phi_2(x_1,y_1,\xi_1)=(x_1-y_1)\xi_1-\frac{(x_1-y_1)^2}{2(T-y_1)}\xi_1=(x_1-y_1)\omega_1\,\,\mbox{\rm with}\,\,\omega_1=(1-\frac{x_1-y_1}{2(T-y_1)})\xi_1.
\end{array}
\end{equation}
Since the Jacobian of the transformation
$\omega_1\mapsto\xi_1$ is $(1-\frac{x_1-y_1}{2(T-y_1)})^{-1}$, we have
\begin{equation}\label{7.6}
\begin{array}{l}
 \int e^{i\Phi_2(x_1,y_1,\xi_1)}\beta_1(\xi) e^{\Phi_3(x_1,y_1,\xi')} \tilde{d}\xi_1\\
=\theta_\lambda(\frac{\xi'}{kT})\int e^{i(x_1-y_1,\omega_1)}\eta\big((2\gamma kT)^{-1}(1-\frac{x_1-y_1}{2(T-y_1)})^{-1}\omega_1\big)(1-\frac{x_1-y_1}{2(T-y_1)})^{-1} e^{\Phi_3(x_1,y_1,\xi')} \tilde{d}\omega_1.
\end{array}
\end{equation}

This is the Schwartz kernel of the pseudo-differential operator over  $(0,T/2)$ with double symbol
\begin{equation}\label{7.7}
\begin{array}{l}
Z:=\theta_\lambda(\frac{\xi'}{kT})\eta\big((2\gamma kT)^{-1}(1-\frac{x_1-y_1}{2(T-y_1)})^{-1}\omega_1\big)(1-\frac{x_1-y_1}{2(T-y_1)})^{-1} e^{\Phi_3(x_1,y_1,\xi')}
\end{array}
\end{equation}
and its simplified symbol $Z_L$ has an asymptotic expansion
\begin{equation}\label{7.8}
\begin{array}{l}
Z_L(x,\omega_1,\xi') \sim\\
\quad\sum_\alpha(\alpha!)^{-1}\partial_{\omega_1}^{\alpha_1} D_{y_1}^{\alpha_1}\{\eta((2\gamma kT)^{-1}(1-\frac{x_1-y_1}{2(T-y_1)})^{-1}\omega_1)(1-\frac{x_1-y_1}{2(T-y_1)})^{-1} e^{\Phi_3(x_1,y_1,\xi')}\} \mid_{y_1=x_1}.
\end{array}
\end{equation}
Here note that $\partial_{\omega_1}$ gives an effect of amplifying by $(kT)^{-1}$ and that of $D_{y_1}$ is amplifying by $\mbox{max}\{k^{1/2},T^{-1}\}$. These total effect can be transformed into amplifying by $O((k^{1/2}T)^{-1})$ if we take
$k$ and $T$ to satisfy $k\ge T^{-2}$.

Thus, we have
\begin{equation}\label{7.9}
\mathcal{P}_k^* \mathcal{S}_k'=I+R_b,
\end{equation}
where $R_b$ is a bounded linear operator on $L^2((0,T/2)\times{\Bbb R}^{n-1})$ with its operator norm estimated by $O(k^{-1/6})$ if we take $k$ and $T$ to satisfy
$k\ge T^{-3}$. Therefore we can modify
$\mathcal{S}_k'$ to obtain a fundamental solution of $\mathcal{P}_k^*$ which can have the same estimate as that of $\mathcal{S}_k'$ if we can show that $\mathcal{S}_k'$ is a bounded linear operator on $L^2((0,T/2)\times{\Bbb R}^{n-1})$. Also, this relation on the estimates of $\mathcal{S}_k'$ and its fundamental solution is true for their derivatives.

\section{Estimate of $S_k'(x_1,y_1;\xi')$}
Let $x_0\in U$ be such that $\dot{B}(\xi')=B(x_0,\xi')$ is diagonalizable for any $\xi'\in{\Bbb R}^{n-1}\setminus0$. Put $\dot{B}_b=\dot{B}+ik(y_1-T)$. Then
\begin{equation*}
\begin{array}{l}
M(y_1,\xi')^{-1}=(T-y_1)(\xi_1-\dot{B}_b)^{-1},\,\,\,
\Phi_1(x_1,y_1,\xi')=-\frac{(x_1-y_1)^2}{2(T-y_1)}
(\xi_1-\dot{B}).
\end{array}
\end{equation*}
By the {\sl further assumptions} and spectral decomposition of $\dot{B}$, we have
\begin{equation}\label{6.3}
\begin{array}{l}
e^{i\Phi}M(y_1,\xi')^{-1}\\
=(T-y_1)e^{ i[(x_1-y_1)-\frac{(x_1-y_1)^2}{2(T-y_1)}]\xi_1}e^{i\frac{(x_1-y_1)^2}{2(T-y_1)}(\dot{B}_R+i\dot{B}_I)}\cdot\sum_{h=1}^s\{\varphi(\lambda_h)(\xi_1-\tilde{\lambda}_h)^{-1}P_h \},
\end{array}
\end{equation}
where $\tilde{\lambda}_h=\lambda_h+ik(y_1-T)$, $P_h=P_h(x_0,\xi')$ is the projection defined before by (\ref{projection}) and $\varphi(\lambda)$ is defined by
\begin{equation}
\varphi(\lambda)=\exp[i(x_1-y_1)^2\lambda/(2(T-y_1))].
\end{equation}

Let $\lambda_j=\lambda(x_0,\xi')=\nu_j+i\mu_j$ with $\nu_j=\nu_j(x_0,\xi'),\,\mu_j=\mu_j(x_0,\xi')\in{\Bbb R}$ for $1\leq j\leq s$ be the eigenvalues of $\dot{B}(\xi')=\dot{B}(x_0,\xi')$ with multiplicity $m_j$ satisfying $m_1+\cdots+m_s=N$. For the estimate of \eqref{6.3}, we need to estimate
\begin{equation*}
\begin{array}{l}
\int \eta(\frac{\xi_1}{2\gamma kT}) e^{i\{(x_1-y_1)-\frac{(x_1-y_1)^2}{2(T-y_1)}\}\xi_1}e^{\frac{(x_1-y_1)^2}{2(T-y_1)}(i\nu_j-\mu_j)}\{ \xi_1-\nu_j-i(\mu_j-k(T-y_1))\}^{-1}\tilde{d}\xi_1.
\end{array}
\end{equation*}
By taking $\gamma$ large to satisfy $|\nu_j|\le1/(4\gamma kT)\,\,(1\le j\le N)$ so that each of $|(\xi_1+\nu_j)/(2\gamma kT)|\le 1/2\,\,(1\le j\le N)$ implies $|\xi_1/(2\gamma kT)|\le 1$, it is enough to estimate
\begin{equation*}
\begin{array}{l}
n_j(x_1,y_1;x_0,\xi'):=\\
\int \eta(\frac{\xi_1}{2\gamma kT}) e^{i\{(x_1-y_1)-\frac{(x_1-y_1)^2}{2(T-y_1)}\}\xi_1}e^{-\frac{(x_1-y_1)^2}{2(T-y_1)}\mu_j}\{ \xi_1-i(\mu_j-k(T-y_1))\}^{-1}\tilde{d}\xi_1
\end{array}
\end{equation*}
for each $1\leq j\leq N$.
For the estimate of each $n_j(x_1,y_1,\xi')$ observe that
\begin{equation}\label{6.4}
\begin{array}{l}
|\int_{|\xi_1|\leq 2\gamma kT}(1- \eta(\frac{\xi_1}{2\gamma kT})) e^{i\{(x_1-y_1)-\frac{(x_1-y_1)^2}{2(T-y_1)}\}\xi_1} e^{-\frac{(x_1-y_1)^2}{2(T-y_1)}\mu_j}\{ \xi_1-i(\mu_j-k(T-y_1))\}^{-1}\tilde{d}\xi_1|\\
\leq (2\pi)^{-1}\int_{\gamma kT\le|\xi_1|\leq 2\gamma kT}\frac{1}{|\xi_1|}d\xi_1\\
=\pi^{-1}
\log 2.
\end{array}
\end{equation}
So we only need to estimate
$$\int_{-2\gamma kT}^{2\gamma kT} e^{i\{(x_1-y_1)-\frac{(x_1-y_1)^2}{2(T-y_1)}\}\xi_1} e^{-\frac{(x_1-y_1)^2}{2(T-y_1)}\mu_j(\xi')}\{ \xi_1-i(\mu_j-k(T-y_1))\}^{-1}\tilde{d}\xi_1.$$
To begin estimating this, set
\begin{equation}\label{6.5}
\begin{array}{l}
I=\int_{-2\gamma kT}^{2\gamma kT} e^{i\{(x_1-y_1)-\frac{(x_1-y_1)^2}{2(T-y_1)}\}\xi_1}e^{-\frac{(x_1-y_1)^2}{2(T-y_1)}\mu_j(\xi')}\{ \xi_1-i(\mu_j-k(T-y_1))\}^{-1}\tilde{d}\xi_1
\end{array}
\end{equation}
and $\mu=\mu_j(y_1,\xi')$, $x=x_1-y_1$, $m=T-y_1$, $\zeta=\xi_1$, $\tilde{d}\zeta=\tilde{d}\xi_1$, then
\begin{equation}\label{6.6}
\begin{array}{l}
I=\int_{-2\gamma kT}^{2\gamma kT} e^{i(x-x^2m^{-1}/2)\zeta}e^{-x^2m^{-1}\mu/2}\cdot\{\zeta-i(\mu-km)\}^{-1}\tilde{d}\zeta\\
=e^{-x^2m^{-1}\mu/2}\int_{-2\gamma kT}^{2\gamma kT} e^{i(x-x^2m^{-1}/2)\zeta} \frac{\zeta+i(\mu-km)}{\zeta^2+(\mu-km)^2}\tilde{d}\zeta.
\end{array}
\end{equation}
Now, we consider two cases. For the case $\mu-km\neq 0$.
Observe that
\begin{equation}\label{6.7}
\begin{array}{l}
|\mu-km|\int_{-2\gamma kT}^{2\gamma kT} \frac{1}{\zeta^2+(\mu-km)^2}d\zeta=\tan^{-1}(\frac{\zeta}{|\mu-km|})]_{-2\gamma kT}^{2\gamma kT}\leq\pi.
\end{array}
\end{equation}
So, we only need to estimate
\begin{equation}\label{6.8}
\begin{array}{l}
\frac{J}{2i}:=\frac{1}{2i}\int_{0}^{2\gamma kT} \frac{\zeta}{\zeta^2+(\mu-km)^2}\,e^{i(x-x^2 m^{-1}/2)\zeta}\,d\zeta+\frac{1}{2i}\int_{-2\gamma kT}^{0} \frac{\zeta}{\zeta^2+(\mu-km)^2}\,e^{i(x-x^2 m^{-1}/2)\zeta}\,d\zeta\\
=\int_{0}^{2\gamma kT} \frac{\zeta\cdot\sin[(x-x^2m^{-1}/2)\zeta]}{\zeta^2+(\mu-km)^2}d\zeta\\
=\int_{0}^{2\gamma kT} \frac{\zeta\cdot\sin[b\zeta]}{\zeta^2+a^2}\,d\zeta\\
=\int_{0}^{2\gamma bkT} \frac{\zeta\cdot\sin\zeta}{\zeta^2+b^2a^2}d\zeta,
\end{array}
\end{equation}
where $a=\mu-km$, $b=x-x^2m^{-1}/2$.

If $b=0$, then we just have $J/(2i)=0$. Hence we assume $b\not= 0$. Let $l\in {\Bbb N}\cup\{0\}$ such that $2\gamma bKT-2\pi l<2\pi$.
For any natural number $j\,(j<l)$, we have
\begin{equation}\label{6.9}
\begin{array}{l}
|\int_{2\pi j}^{2\pi j+2\pi} \frac{\zeta\cdot\sin\zeta}{\zeta^2+b^2a^2}d\zeta|\\
=|\int_{2\pi j}^{2\pi j+\pi} \frac{\zeta\cdot\sin\zeta}{\zeta^2+b^2a^2}d\zeta+\int_{2\pi j+\pi}^{2\pi j+2\pi} \frac{\zeta\cdot\sin\zeta}{\zeta^2+b^2a^2}d\zeta|\\
=|\int_{2\pi j}^{2\pi j+\pi} \frac{\zeta\cdot\sin\zeta}{\zeta^2+b^2a^2}d\zeta+\int_{2\pi j}^{2\pi j+\pi} \frac{(\zeta+\pi)\cdot\sin(\zeta+\pi)}{(\zeta+\pi)^2+b^2a^2}d\zeta|\\
=|\int_{2\pi j}^{2\pi j+\pi} \frac{\zeta\cdot\sin\zeta}{\zeta^2+b^2a^2}d\zeta-\int_{2\pi j}^{2\pi j+\pi} \frac{(\zeta+\pi)\cdot\sin\zeta}{(\zeta+\pi)^2+b^2a^2}d\zeta|\\
=|\int_{2\pi j}^{2\pi j+\pi} [\frac{\zeta}{\zeta^2+b^2a^2}-\frac{(\zeta+\pi)}{(\zeta+\pi)^2+b^2a^2}]\cdot\sin\zeta d\zeta\\
=|\int_{2\pi j}^{2\pi j+\pi} \frac{\zeta^2\pi+\zeta\pi^2-b^2a^2\pi}{(\zeta^2+b^2a^2)[(\zeta+\pi)^2+b^2a^2]}\cdot\sin\zeta d\zeta\\
\leq \int_{2\pi j}^{2\pi j+\pi} \frac{2\pi}{(\zeta+\pi)^2} d\zeta.
\end{array}
\end{equation}
Now, we estimate $\frac{J}{2i}$. We have that
\begin{equation}\label{6.10}
\begin{array}{l}
|\frac{J}{2i}|=|\Sigma_{j=1}^{l-1}\int_{2\pi j}^{2\pi j+2\pi} \frac{\zeta\cdot\sin\zeta}{\zeta^2+b^2a^2}d\zeta+\int_{0}^{2\pi} \frac{\zeta\cdot\sin\eta}{\zeta^2+b^2a^2}d\zeta+\int_{2\pi l}^{2\gamma bkT} \frac{\zeta\cdot\sin\zeta}{\zeta^2+b^2a^2}d\zeta|\\
\leq|\Sigma_{j=1}^{l-1}\int_{2\pi j}^{2\pi j+2\pi} \frac{\zeta\cdot\sin\zeta}{\zeta^2+b^2a^2}d\zeta|+4\pi\\
\leq\int_{0}^{\infty} \frac{2\pi}{(\zeta+\pi)^2} d\zeta+4\pi\\
\leq 2+4\pi.
\end{array}
\end{equation}

On the other hand, for the case $\mu-km=0$, we have
\begin{equation}\label{6.11}
\begin{array}{l}
|I|=|e^{-x^2m^{-1}\mu/2}\,\int_{-kT}^{kT} e^{i(x-x^2m^{-1})\zeta}\,\zeta^{-1}d\zeta|\\
=|2\,e^{-x^2m^{-1}\mu/2}\,\int_{0}^{kT} \frac{\sin[(x-x^2m^{-1})\zeta]}{\zeta}d\zeta|.
\end{array}
\end{equation}
The same argument can give $|I|$ bounded by an explicit positive constant.

Since $(1+k(x_1-y_1)^2)^{-m}\lesssim e^{-\frac{(x_1-y_1)^2}{2(T-y_1)}\mu_j(\xi')}$ for $m\geq 0$, we have finally obtained the estimate for $S_k'(x_1,y_1;\xi')=S_k'(x_1,y_1;x_0,\xi')$
\begin{equation}\label{6.12}
 |S_k'(x_1,y_1;\xi')|\lesssim (1+k(x_1-y_1)^2)^{-m}\,T
\end{equation}
uniformly for any $x_1,\,y_1\in [0,T/2]$, $\xi'\in {\Bbb R}^{n-1}\setminus 0$ and of those $x_0\in U$ with the property stated in the first line of this section. Here $T$ in the above estimate comes from the factor $T-y_1$ in \eqref{6.3}.

\section{Estimate of parametrix $\mathcal{S}_k'$}\label{sec8}
Recall that the Schwartz kernel of $\mathcal{S}_k'$ was given by
\begin{equation*}
\begin{array}{ll}
S_k'(x,y)&=\int_{{\Bbb R}^{n-1}}\theta_\lambda(\frac{\xi'}{kT})e^{i(x'-y')\cdot\xi'}\Big(\int \eta(\frac{\xi_1}{2\gamma kT}) e^{i\Phi(x_1,y_1,\xi')}M^{-1}(y_1,\xi')\tilde{d}\xi_1\Big)\tilde{d}\xi'\\
&=\int_{{\Bbb R}^{n-1}}e^{i(x'-y')\cdot\xi'}S_k(x_1,y_1;\xi')\tilde{d}\xi'.
\end{array}
\end{equation*}
In this section, we will estimate $\partial_x^{\alpha}S_k' v$ for $|\alpha|\leq 1$ and $g\in C_0^\infty((0,T/2)\times{\Bbb R}^{n-1})$.
To begin with,
\begin{equation}\label{8.1}
\begin{array}{l}
\|\mathcal{S}_k' g\|_{L^2(\mathbf{R}^{n-1},dx')}=(\int|\mathcal{S}_k'g(x)|^2 dx')^{1/2}\\
=(\int|\iint S_k'(x,y)g(y)dy'dy_1|^2 dx')^{1/2}\\
\leq\int(\int|\int S_k'(x,y)g(y)dy'|^2 dx')^{1/2}dy_1,
\end{array}
\end{equation}
where we used the Minkowski integral inequality in the last inequality.

Now we estimate the inner integral
\begin{equation}\label{8.2}
\begin{array}{l}
\int S_k'(x,y)g(y)dy'\\
=\iint e^{i(x'-y')\xi'}S_k'(x_1,y_1,\xi')g(y)dy'd\xi'\\
=\int e^{ix'\xi'}S_k'(x_1,y_1,\xi')\hat{g}(y_1,\xi')d\xi'\\
=\widetilde{S_k'\hat{g}}(y_1,x),
\end{array}
\end{equation}
where "\,$\widetilde{\cdot}$\," denotes the inverse Fourier transform of "$\cdot\,$".
By the Plancherel theorem and \eqref{6.12}, we can estimate
\begin{equation}\label{8.3}
\begin{array}{l}
(\int|\int S_k'(x,y)g(y)dy'|^2 dx')^{1/2}\\
=(\int|S_k'(x_1,y_1,\xi')\hat{g}(y_1,\xi')|^2 d\xi')^{1/2}\\
\leq  T [1+k(x_1-y_1)^2]^{-1}\cdot(\int|g(y)|^2 dy')^{1/2}.
\end{array}
\end{equation}
So, we have
\begin{equation}\label{8.4}
\begin{array}{l}
\|\mathcal{S}_k' g\|^2_{L^2((0,T/2)\times{\Bbb R}^{n-1})}=\int_0^T\|\mathcal{S}_k' g(x)\|^2_{L^2(\mathbf{R}^{n-1},dx')}dx_1\\
\leq \int \big[\int T [1+k(x_1-y_1)^2]^{-1}\cdot(\int|g(y)|^2 dy')^{1/2}\, dy_1\big]^{2}dx_1\\
\leq T^2(\int[1+k x_1^2]^{-1}dx_1)^2\|g\|^2_{L^2(\mathbf{R}^{n})}\\
\leq \frac{T^2}{k}\|g\|^2_{L^2((0,T/2)\times{\Bbb R}^{n-1})},
\end{array}
\end{equation}
where we used the Young inequality in the second inequality.

Let $\mathcal{S}_{k,1}'$ and $\mathcal{S}_{k,0}'$ be the operators with Schwartz's kernels $S_{k,1}'(x,y)$ and $S_{k,0}'(x,y)$, respectively. Then, the same argument gives that for $|\alpha|\leq 1$
\begin{equation}\label{8.5}
\begin{array}{l}
\|\partial_x^{\alpha}\mathcal{S}_{k,1}' g\|_{L^2((0,T/2)\times{\Bbb R}^{n-1})}\leq \frac{T}{\sqrt{k}}(kT)^{|\alpha|}\|g\|_{L^2((0,T/2)\times{\Bbb R}^{n-1})}.
\end{array}
\end{equation}
For the estimate of $\partial_x^{\alpha}\mathcal{S}_{k,0}'v$ for $|\alpha|\leq 1$, we surely have from \eqref{5.14} that
\begin{equation}\label{8.6}
\begin{array}{l}
\|\partial_x^{\alpha}\mathcal{S}_{k,0}'g\|_{L^2((0,T/2)\times{\Bbb R}^{n-1})}\leq \frac{T}{\sqrt{k}}(kT)^{|\alpha|}\|g\|_{L^2((0,T/2)\times{\Bbb R}^{n-1})}.
\end{array}
\end{equation}
Thus for $|\alpha|\le 1$, we have
\begin{equation}\label{8.7}
\begin{array}{l}
\|\partial_x^{\alpha}\mathcal{S}_{k}'g\|_{L^2((0,T/2)\times{\Bbb R}^{n-1})}\leq \frac{T}{\sqrt{k}}(kT)^{|\alpha|}\|g\|_{L^2((0,T/2)\times{\Bbb R}^{n-1})}.
\end{array}
\end{equation}
Combining \eqref{8.7} with what was mentioned in the last paragraph of Section 6, we obtain for $v\in C_0^\infty((0,T/2)\times{\Bbb R}^{n-1})$ that
\begin{equation}\label{8.8}
\begin{array}{rl}
\sum_{|\beta|\le 1}T^{-\frac{1}{2}}(kT)^{\frac{1}{2}-|\beta|}\Vert \partial_{x}^\beta v\Vert_{L^2((0,T/2)\times{\Bbb R}^{n-1})}&\lesssim T\|\mathcal{P}_k^*v\|_{L^2((0,T/2)\times{\Bbb R}^{n-1})}\\
&\lesssim \|P_{k,b}^*v\|_{L^2((0,T/2)\times{\Bbb R}^{n-1})}.
\end{array}
\end{equation}
Since $S_k'\partial_x^\alpha$ can have a similar estimate as \eqref{8.7} and its adjoint is $(-1)^{|\alpha|}\partial_x^\alpha S_k^*$, we can see likewise (\ref{8.8})
\begin{equation}\label{8.9}
\begin{array}{l}
\sum_{|\beta|\le 1}T^{-\frac{1}{2}}(kT)^{\frac{1}{2}-|\beta|}\
\Vert \partial_{x}^\beta v\Vert_{L^2((0,T/2)\times{\Bbb R}^{n-1})}\lesssim \|P_{k,b}v\|_{L^2((0,T/2)\times{\Bbb R}^{n-1})}.
\end{array}
\end{equation}

\section{Estimate of $A_k$}\label{sec9}
First of all since $P_{k,g}$ is an elliptic pseudo-differential operator of order 1 with large parameter $k$, we can easily have
\begin{equation}\label{9.4}
\begin{array}{l}
\|\partial_x^{\alpha}v\|_{L^2((0,T/2)\times{\Bbb R}^{n-1})}\lesssim (kT)^{-1+|\alpha|-\ell}\sum_{|\beta|\le\ell}\|\partial_x^\beta P_{k,g}v\|_{L^2((0,T/2)\times{\Bbb R}^{n-1})}\,\,\,
(\ell=0,1,\,|\alpha|\leq 2)
\end{array}
\end{equation}
for any $v\in C_0^\infty((0,T/2)\times{\Bbb R}^{n-1})$ by the standard argument (see for instance Chapter 2 of \cite{shubin}).

Combining \eqref{8.9} and \eqref{9.4}, we have for any $v\in C_0^\infty((0,\frac{T}{2})\times{\Bbb R}^{n-1})$ that
\begin{equation}\label{9.6}
\begin{array}{l}
\sum_{\alpha+|\beta|\le 2}(kT)^{3-2\alpha-2|\beta|}T^{-1}\Vert \partial_{x_1}^\alpha\partial_{x'}^\beta v\Vert^2
\lesssim \Vert  A_kv\Vert^2
\end{array}
\end{equation}
and hence we have the following Carleman estimates for constant coefficients partial differential operator $\dot{\mathcal{L}}$.

\begin{theorem}\label{thm9.1}
There exist positive constants $T_0$, $k_0$ and $c$ such
that for $0<T\leq T_0$, we have for $k\geq k_0$ that
\begin{equation}\label{9.7}
\begin{array}{l}
(kT^2)^{-1}\int_0^Te^{k(x_1-T)^2}\|\partial^{2}v\|^2dx_1+k\int_0^Te^{k(x_1-T)^2}\|\partial v\|^2dx_1+k^3T^2\int_0^Te^{k(x_1-T)^2}\|v\|^2dx_1\\
\leq c\int_0^Te^{k(x_1-T)^2}\|\dot{\mathcal{L}}v\|^2dx_1
\end{array}
\end{equation}
for any $v\in C_0^\infty((0,T/2)\times{\Bbb R}^{n-1})$, where $\|\cdot\|^2=(\cdot,\cdot)$
is the $L^2({\Bbb R}^{n-1})$ norm.
\end{theorem}

\section{Carleman estimates by partition of unity}\label{sec10}
We will use a partition of unity and Theorem \ref{thm9.1} to prove Theorem \ref{thm2.2}. To begin with we first note that we are back to use the notations used to describe (\ref{tildeL}). let $\vartheta_0\in C^\infty_0({\Bbb R})$ and $0\leq\vartheta\leq 1$ such that
\begin{equation}\label{10.1}
\vartheta_0(t)=
\begin{cases}
\begin{array}{l}
1\,\,\,(|t|\leq 1)\\
0\,\,\,(|t|\geq 3/2).
\end{array}
\end{cases}
\end{equation}
Also let  $\vartheta(\tilde{x})=\vartheta_0(\tilde{x}_1)\cdots\vartheta_0(\tilde{x}_n)$. Then, we have
\begin{equation}\label{10.2}
\vartheta(\tilde{x})=
\begin{cases}
\begin{array}{l}
1\,\,\,(\tilde{x}\in Q_1(0))\\
0\,\,\,(\tilde{x}\in{\Bbb R}^N\setminus Q_{3/2}(0)),
\end{array}
\end{cases}
\end{equation}
where $Q_{r}=Q_{r}(0)=\{x\in{\Bbb R}^n\,:\,|x_j|\leq r,\,j=1,2,\cdots ,n\}$.
Further for $\mu\geq 1$ and $g\in {\Bbb Z}^n$, define
$$\tilde{x}_g=g/\mu$$
and set
$$\vartheta_{g,\mu}(\tilde{x})=\vartheta(\mu \tilde{x}-g).$$
Then we have
\begin{equation}\label{10.3}
\begin{array}{l}
{\rm supp}\,\vartheta_{g,\mu}\subset Q_{3/2\mu}(x_g)\subset Q_{2/\mu}(x_g)
\end{array}
\end{equation}
and
\begin{equation}\label{10.4}
\begin{array}{l}
|D^k\vartheta_{g,\mu}|\leq c_1\mu^k(\chi_{ Q_{3/2\mu}(\tilde{x}_g)}-\chi_{Q_{1/\mu}(\tilde{x}_g)})\,\,\,(k=0,1,2),
\end{array}
\end{equation}
where $c_1\geq 1$ depends only on $n$.

For $g\in {\Bbb Z}^n$, let
$$\mathcal{A}_g=\{g'\in {\Bbb Z}^n\,:\,{\rm supp}\vartheta_{g',\mu}\cap{\rm supp}\vartheta_{g,\mu}\neq\varnothing\},$$
then the number of $\mathcal{A}_g$ will only depend on $n$ but  not on $\mu$. Thus, we can define
\begin{equation}\label{10.5}
\begin{array}{l}
\bar{\vartheta}_{\mu}(\tilde{x}):=\sum_{g\in {\Bbb Z}^N}\vartheta_{g,\mu}(\tilde{x})\geq 1\,\,\,(\tilde{x}\in {\Bbb R}^n).
\end{array}
\end{equation}
From \eqref{10.4}, we have
\begin{equation}\label{10.6}
\begin{array}{l}
|D^k\bar{\vartheta}_{\mu}|\leq c_2\mu^k,
\end{array}
\end{equation}
where $D^k$ denotes all the $k$ th order derivatives and $c_2\geq 1$ depends only on $n$.

Define
$$\eta_{g,\mu}(\tilde{x})=\vartheta_{g,\mu}(\tilde{x})/\bar{\vartheta}_{\mu}(\tilde{x})\,\,\,(\tilde{x}\in {\Bbb R}^n).$$
Then, we have
\begin{equation}\label{10.7}
\begin{cases}
\begin{array}{l}
\sum_{g\in {\Bbb Z}^N}\eta_{g,\mu}(\tilde{x})= 1\,\,\,(\tilde{x}\in {\Bbb R}^N)\\
{\rm supp}\,\eta_{g,\mu}\subset Q_{3/2\mu}(\tilde{x}_g)\subset Q_{2/\mu}(\tilde{x}_g)\\
|D^k\eta_{g,\mu}|\leq c_3\mu^k\chi_{ Q_{3/2\mu}(\tilde{x}_g)}\,\,\,(k=0,1,2),
\end{array}
\end{cases}
\end{equation}
where $c_3\geq 1$ depends only on $n$.

\medskip

Now define $\widetilde{L}_g(\tilde{x},\partial)u=\sum_{1\leq j,\,\ell\leq n}\sum_{1\leq\beta\leq N}\widetilde{C}_{\alpha\beta}^{pq}(\tilde{x}_g)\partial_p\partial_qu_{\beta}$. Then from Theorem \ref{thm9.1} we have the following proposition.

\begin{pr}\label{thm10.1}
There exist positive constants $T_0$, $k_0$, $r$ and $c$ such
that for $0<T\leq T_0$, we have for $k\geq k_0$ that
\begin{equation}\label{10.8}
\begin{array}{l}
(kT^2)^{-1}\int_0^Te^{k(\tilde{x}_1-T)^2}\|\partial^{2}v\|^2d\tilde{x}_1+k\int_0^Te^{k(\tilde{x}_1-T)^2}\|\partial v\|^2d\tilde{x}_1+k^3T^2\int_0^Te^{k(\tilde{x}_1-T)^2}\|v\|^2d\tilde{x}_1\\
\leq c\int_0^Te^{k(\tilde{x}_1-T)^2}\|\widetilde{L}_g v\|^2d\tilde{x}_1
\end{array}
\end{equation}
for all $v(\tilde{x}_1,\tilde{x}')\in C_0^{\infty}({\Bbb R}^n)$ with $\mbox{\rm
supp}\,v\subset\mathcal{U}:=\{(\tilde{x}_1,\tilde{x}')\,:\,\tilde{x}_1\in(0,T/2),\,|\tilde{x}'|<r\}\subset\widetilde{U}$.
\end{pr}
Now let ${\rm supp}(v)\subset \mathcal{U}/2$ and apply \eqref{10.8} to $v\eta_{g,\mu}$. Then we have
\begin{equation}\label{10.9}
\begin{array}{rl}
&(kT^2)^{-1}\int_0^Te^{k(\tilde{x}_1-T)^2}\|D^{2}(v\eta_{g,\mu})\|^2d\tilde{x}_1+k\int_0^Te^{k(\tilde{x}_1-T)^2}\|D(v\eta_{g,\mu})\|^2d\tilde{x}_1\\
&+k^3T^2\int_0^Te^{k(\tilde{x}_1-T)^2}\|v\eta_{g,\mu}\|^2d\tilde{x}_1\\
\leq & c\int_0^Te^{k(\tilde{x}_1-T)^2}\|\widetilde{L}_g(v\eta_{g,\mu})\|^2d\tilde{x}_1.
\end{array}
\end{equation}
Note here we have
\begin{equation}\label{10.10}
\begin{array}{rl}
|\widetilde{L}_{g}(\tilde{x},\partial)(v\eta_{g,\mu})|\lesssim &|\widetilde{L}(\tilde{x},\partial)(v\eta_{g,\mu})|\\
&+|\widetilde{L}(\tilde{x},\partial)(v\eta_{g,\mu})-\widetilde{L}_{g}(\tilde{x},\partial)(v\eta_{g,\mu})|\\
\lesssim &\eta_{g,\mu}|\widetilde{L}(\tilde{x},\partial)(v)|+\mu^{-1}|D^2v|\,\chi_{Q_{\frac{2}{\mu}}(\tilde{x}_g)}\\
&+\mu|Dv|\,\chi_{Q_{\frac{2}{\mu}}(\tilde{x}_g)}+\mu^2|v|\,\chi_{Q_{\frac{2}{\mu}}(\tilde{x}_g)}.
\end{array}
\end{equation}
By summing up all (\ref{10.10}) with respect to $g\in {\Bbb Z}^n$ and using from \eqref{10.7}, \eqref{10.9},\eqref{10.10}, we have
\begin{equation}\label{10.11}
\begin{array}{rl}
&[(kT^2)^{-1}-c_1\mu^{-2}]\int_0^Te^{k(\tilde{x}_1-T)^2}\|D^{2}v\|^2d\tilde{x}_1+[k-c_2\mu^2]\int_0^Te^{k(\tilde{x}_1-T)^2}\|Dv\|^2d\tilde{x}_1\\
&+[k^3T^2-c_3k\mu^2-c_4\mu^4]\int_0^Te^{k(\tilde{x}_1-T)^2}\|v\|^2d\tilde{x}_1\\
\lesssim & c\int_0^Te^{k(\tilde{x}_1-T)^2}\|\widetilde{L}v|^2d\tilde{x}_1.
\end{array}
\end{equation}
Taking $c_2\mu^2=k/2$, $4T=c_2/c_1$ and $k$ large, we can derive that
\begin{equation}\label{10.12}
\begin{array}{rl}
&(kT^2)^{-1}\int_0^Te^{k(\tilde{x}_1-T)^2}\|D^{2}v\|^2d\tilde{x}_1+k\int_0^Te^{k(\tilde{x}_1-T)^2}\|Dv\|^2d\tilde{x}_1\\
&+k^3T^2\int_0^Te^{k(\tilde{x}_1-T)^2}\|v\|^2d\tilde{x}_1\\
\lesssim & c\int_0^Te^{k(\tilde{x}_1-T)^2}\|\widetilde{L}v\|^2d\tilde{x}_1.
\end{array}
\end{equation}

\end{document}